\documentclass[10pt,conference,onecolumn,letterpaper]{IEEEtran}
\usepackage{epsfig,amsmath,hhline}

\usepackage{times}

\usepackage{latexsym}
\usepackage{amsmath}
\usepackage{nccbbb}
\usepackage{amsthm}
\usepackage{url}
\usepackage{algorithm}
\usepackage{algorithmic}
\usepackage{dsfont}

\newtheorem{theorem}{Theorem}
\newtheorem{fact}{Fact}
\newtheorem{lemma}{Lemma}

\theoremstyle{definition}
\newtheorem{defn}{Definition}
\newcommand{\ie}{\textsl{i.e.} }

\newcommand{\Var}{\mathrm{Var}}

\usepackage{graphicx}
\usepackage{amssymb}
\usepackage{epstopdf}
\title{Heavy Traffic Queue Length Behavior in Switches with Reconfiguration Delay}

\author{
\IEEEauthorblockN{Chang-Heng Wang}\\
\IEEEauthorblockA{
University of California, San Diego\\
Email: chw009@ucsd.edu}
\and
\IEEEauthorblockN{Siva Theja Maguluri}\\
\IEEEauthorblockA{
Georgia Institute of Technology\\
Email: siva.theja@gatech.edu}
\and
\IEEEauthorblockN{Tara Javidi}\\
\IEEEauthorblockA{
University of California, San Diego\\
Email: tjavidi@ucsd.edu}
}

\begin{document}

\maketitle

\begin{abstract}

Optical switches have been drawing attention due to their large data bandwidth and low power consumption. However, scheduling policies need to account for the schedule reconfiguration delay of optical switches to achieve good performance. The Adaptive MaxWeight policy achieves optimal throughput for switches with nonzero reconfiguration delay, and has been shown in simulation to have good delay performance. In this paper, we analyze the queue length behavior of a switch with nonzero reconfiguration delay operating under the Adaptive MaxWeight. We first show that the Adaptive MaxWeight policy exhibits a weak state space collapse behavior in steady-state, which could be viewed as an inheritance of the MaxWeight policy in a switch with zero reconfiguration delay. We then use the weak state space collapse result to obtain a steady state delay bound under the Adaptive MaxWeight algorithm in heavy traffic by applying a recently developed drift technique. The resulting delay bound is dependent on the expected schedule duration. We then derive the relation between the expected schedule duration and the steady state queue length through drift analysis, and obtain asymptotically tight queue length bounds in the heavy traffic regime.

\end{abstract}

\section{Introduction}

Modern data centers aggregate huge amount of computing and storage resource to support high demanding applications such as cloud computing, large-scale web applications, big data analytics, etc. With the ever increasing number of resources and the communication demand between these resources, the interconnecting networks face stringent performance challenge. Optical switches emerge as a promising candidate to address this challenge since they can support higher data bandwidth relatively easier than traditional electronic switches, and also have lower power consumption. However, optical switches pose another challenge different from traditional electronic switches that makes it difficult to directly substitute electronic switches: optical switches typically exhibit a delay following each schedule reconfiguration, during which no packet transmission could occur~\cite{nistica, ROADM}. This delay is referred as the \textbf{reconfiguration delay}, and it makes the switch scheduling problem more difficult. For example, it is known that with the reconfiguration delay, the well-known MaxWeight policy~\cite{Tassiulas, throughput100} is not even throughput optimal.



For the scheduling of optical switches, or in general, switches with reconfiguration delay, many works in the literature consider the problem in ``quasi-static'' sense: decomposing the traffic demand into a sequence of efficient schedules (in the sense of minimizing service time and number of schedule reconfigurations) for a predetermined time horizon, e.g. \cite{ADJUST, submodular, Solstice, reactor}. The performance of such solutions is usually limited by the duration of the time horizon, and may require some prior knowledge of the traffic arrival rate to achieve good performance. In contrast, \cite{AMW_Infocom, MWMH, SCB} pose the dynamic scheduling problem for switches with reconfiguration delay. The proposed Adaptive MaxWeight~\cite{AMW_Infocom} (or a similar variant, Switching Curve Based policy~\cite{SCB}) makes schedule decision at every time slot, as opposed to scheduling over a time horizon. The key idea of the Adaptive MaxWeight is to reconfigure the schedule only when the current schedule is not good enough, which would be described with more detail in Section~\ref{Sec:AMW}. The idea has also been generalized in~\cite{AMW_arxiv} to introduce a large class of scheduling policies for switches with reconfiguration delay. These policies have been shown to guarantee throughput optimality under mild assumptions on arrival traffic, which means that the policies could guarantee finite expected queue lengths whenever any other policy could achieve that.

In this work, we focus on the analysis of the Adaptive MaxWeight policy. Beyond throughput optimality, it is desirable to further consider the delay behavior (or equivalently, queue length behavior) in order to better evaluate the performance. However, similar to the MaxWeight policy in switches without reconfiguration delay, it is unclear whether an exact expression for the steady state queue length behavior could be obtained. Therefore one usually approaches the performance evaluation by increasing the number of queues in the system, or increasing the arrival rate close to the boundary of the capacity region to study the queue length scaling behavior. In this paper, we consider the queue length behavior of switches with reconfiguration delay (operated under the Adaptive MaxWeight) in the heavy traffic regime. The arrival traffic considered here approaches to a limit where all input ports and all output ports are saturated. 
Under this condition, the steady state queue length under the Adaptive MaxWeight exhibits a behavior similar to the state space collapse (SSC) as introduced in~\cite{MW_scaling}, but in a weaker sense. With the SSC result, we then utilize the drift technique introduced in~\cite{Atilla_drift}, and the Lyapunov function proposed in~\cite{MW_scaling} to obtain a steady state queue length bound. The challenge for the Adaptive MaxWeight case, compared to the MaxWeight under zero reconfiguration delay, is that the steady state queue length bound thus obtained is dependent on the expected schedule duration. One of the main contributions of this paper is in identifying the relation between the expected schedule duration and the expected queue length in steady state, and apply this relation to obtain asymptotically tight bounds for the expected steady state queue length under Adaptive MaxWeight.

The rest of the paper is organized as follows. The switch with reconfiguration delay model and the notion of throughput optimality is introduced in Section~\ref{Sec:model}. In Section~\ref{Sec:AMW}, we briefly introduce the Adaptive MaxWeight policy and some of its properties. We then present our main result regarding heavy traffic queue length behavior in Section~\ref{Sec:heavy}. We first establish the multiplicative state space collapse of the Adaptive MaxWeight. With the state space collapse result, we then establish a queue length upper bound for the steady state queue length in the heavy traffic regime, which is dependent on the expected schedule duration. Using the schedule weight as a Lyapunov function, we further derive the relation between the expected schedule duration and the expected steady state queue length through drift analysis, and then derive the scaling of steady state queue length in heavy traffic. Section~\ref{Sec:simulation} presents some simulation results for the Adaptive MaxWeight, in the effort of characterizing the scaling of the expected queue length with respect to some system parameters, and in comparison with the scaling derived in this paper. Finally, we conclude with a summary and some future directions in Section~\ref{Sec:conclusion}.

\section{System Model}\label{Sec:model}
\subsection{Switch Model and Arrival Traffic}

The model considered in this paper is an $n \times n$ input-queued switch, which has $n$ input ports and $n$ output ports. Each input port maintains $n$ separate queues (either physically or virtually), each storing packets destined to one output port. We denote the queue storing packets at input $i$ and destined for output $j$ with the pair $(i, j)$. This model is also known as the input-queued switch.

The system considered is assumed to be time-slotted, with the time indexed as $t \in \bbbn_{+} = \{0, 1, 2, \dots \}$. Each slot duration is the transmission time of a single packet, which is assumed to be a fixed value. Let $a_{ij}(t)$ be the number of packets arrived at queue $(i, j)$ at time $t$. Let $q_{ij}(t)$ be the number of packets in the queue $(i,j)$ at the beginning of the time slot $t$. Write $\mathbf{a}(t) = [a_{ij}(t)], \mathbf{q}(t) = [q_{ij}(t)]$, where $\mathbf{a}(t), \mathbf{q}(t) \in \bbbn_+^{N \times N}$.

We assume the arrival processes $a_{ij}(t)$ to be independent over $i, j \in \{1, 2, \dots, N\}, i\neq j$. 
For each queue $(i,j)$, the arrival process $a_{ij}(t)$ is assumed to be i.i.d. across time slots, with mean $\bbbe[a_{ij}(t)] = \lambda_{ij}$ and variance $\Var(a_{ij}(t)) = \sigma_{ij}^2$. We also assume that $a_{ij}(t)$ has a finite support, \ie $\exists \ a_{\max} < \infty$ such that $a_{ij}(t) \leq a_{\max}$.

\subsection{Schedules and Reconfiguration Delay}

Let $\mathbf{s}(t) \in \{0, 1\}^{N \times N}$ denote the schedule at time slot $t$, which indicates the queues that are being scheduled by the switch. We set $s_{ij}(t) = 1$ if queue $(i,j)$ is scheduled at time $t$, and $s_{ij}(t) = 0$ otherwise.

The feasible schedules for the network are determined by the network topology and physical constraints on simultaneous data transmissions. We let $\mathcal{S} \subset \bbbr^{n^2}$ denote the set of all feasible schedules, \ie $\mathbf{s}(t) \in \mathcal{S}$ for all $t$. We assume at any $t$ each input port can only transmit to at most one output port, and each output port can only receive from at most input port, i.e. $\sum_i s_{ij}(t) \leq 1, \sum_j s_{ij}(t) \leq 1$. This schedule constraint determines the set of feasible schedules $\mathcal{S}$.

Upon reconfiguring a schedule, the network incurs a reconfiguration delay, during which no packet could be transmitted. We make this notion formal through the following two definitions:

\begin{defn}
Let $\{ t_k^S \}_{k=0}^{\infty}$ denote the time instances when the schedule is reconfigured. The schedule between two schedule reconfiguration time instances remains the same, \ie \[
\mathbf{s}(\tau) = \mathbf{s}(t_k^S), \ \ \forall \tau \in [t_k^S, t_{k+1}^S - 1]
\]
\end{defn}

\begin{defn}
Let $\Delta_r$ be the reconfiguration delay associated with reconfiguring the schedule of the network. During the period of schedule reconfiguration, \ie $\forall t \in \cup_{k=0}^{\infty} [t_k^S, t_k^S + \Delta_r]$, the switch does not serve any of its queues. We assume the reconfiguration delay to be an integer multiple of a time slot.
\end{defn}

Let $r(t)$ denote the time for the switch remaining in the reconfiguration delay, with $r(t) = 0$ indicating that the switch is not in reconfiguration at time $t$. Therefore, for $t \in [t_k^S, t_k^S + \Delta_r]$, for some $k \in \bbbn_+$, we have $r(t) = \Delta_r - (t - t_k^S)$; and $r(t) = 0$ for all other $t$.

With the above definitions, we may then write the queue dynamics for any queue $(i,j)$ as
\begin{align}
q_{ij}(t+1) =& \Big[ q_{ij}(t) + a_{ij}(t) - s_{ij}(t)\mathds{1}_{\{r(t) = 0\}}  \Big]^+ \nonumber \\
=& q_{ij}(t) + a_{ij}(t) - s_{ij}(t)\mathds{1}_{\{r(t) = 0\}} + u_{ij}(t)
\end{align}
where $\mathds{1}_{E}$ is the indicator function of event $E$, and $[x]^+ = \max\{x, 0\}$. Note $u_{ij}(t) \in \{0, 1\}$ is the unused service of queue $(i,j)$ when the queue is empty. A useful property of the unused service is that $u_{ij}(t) = 1$ only when $q_{ij}(t+1) = 0$, or in other words, $u_{ij}(t)q_{ij}(t+1) = 0$.

The schedule at each time slot $\mathbf{s}(t)$ is determined by a scheduling policy. In this paper, we consider scheduling policies that determine the schedule at time $t$ based on the queue length $\mathbf{q}(t)$ and the previous schedule $\mathbf{s}(t-1)$. Under this type of policy, the process $\{\mathbf{X}(t)\}_{t=0}^{\infty}$ with $\mathbf{X}(t) = (\mathbf{q}(t), \mathbf{s}(t), r(t)) \in \bbbn_+^{n^2} \times \{0,1\}^{n^2} \times \{0, 1, \dots, \Delta_r\} \triangleq \mathcal{X}$ that describes the switch model is then a discrete time Markov chain.

\subsection{Stability and Capacity Region}

A queue $(i, j)$ is strongly stable if its queue length $q_{ij}(t)$ satisfies:  
\[
\limsup\limits_{t \rightarrow \infty} \frac{1}{t} \sum_{\tau=1}^t\mathbb{E}\{ q_{ij}(\tau) \} < \infty
\]
and we say the system of queues is stable if queue $(i, j)$ is strongly stable for all $i, j \in \{1, 2, \dots, N\}$. A scheduling policy is said to stabilize the system if the system is stable under that scheduling policy for a given traffic rate matrix. With this notion of stability, we define the capacity region $\mathcal{C}$ of the network as the set of all traffic rate matrices such that there exists a scheduling policy which stabilizes the system.

The capacity region is given by the convex hull of the feasible schedules $\mathcal{S}$, that is 
\begin{align*}
\mathcal{C} =& \Big\{ \sum_{\mathbf{s} \in \mathcal{S}} \alpha_{\mathbf{s}} \mathbf{s} : \sum_{\mathbf{s} \in \mathcal{S}} \alpha_{\mathbf{s}} < 1, \ \alpha_{\mathbf{s}} \geq 0, \ \forall \mathbf{s} \in \mathcal{S} \Big\}
\end{align*}
In this paper, we are interested in traffic rate matrices on the outer boundary of $\mathcal{C}$ where all ports are saturated, \ie
\begin{align*}
\mathcal{F} = \Big\{ \boldsymbol{\lambda} \in \bbbr^{n^2} : \sum_{i} \lambda_{ij} = 1, \sum_{j} \lambda_{ij} = 1, \forall i,j \in \{1, 2, \dots, n \} \Big\}
\end{align*}

For any traffic rate matrix $\boldsymbol{\lambda} \in \mathcal{C}$, we say that $\boldsymbol{\lambda}$ is admissible, and define the load of the traffic as $\rho(\boldsymbol{\lambda}) = \max\{r: \boldsymbol{\lambda} \in r\mathcal{C}, \ 0 < r < 1 \}$.

We say that a scheduling policy is throughput optimal if it stabilizes the system for any traffic rate matrix $\boldsymbol{\lambda} \in \mathcal{C}$.

\section{Adaptive MaxWeight} \label{Sec:AMW}

It is known that for switches without reconfiguration delay, the MaxWeight policy is throughput optimal~\cite{throughput100} and has optimal delay scaling in the heavy traffic regime~\cite{MW_scaling, MW_scaling_incompletely_saturated}. However, with the presence of reconfiguration delay, the MaxWeight policy is not even throughput optimal since it does not account for the overhead of frequent schedule reconfiguration.

The Adaptive MaxWeight scheduling policy is presented in Algorithm~\ref{Alg}. The main idea behind the Adaptive MaxWeight is to reconfigure the schedule when the current schedule is not ``good'' enough. Using the schedule weight as the measure of a schedule, the Adaptive MaxWeight computes the schedule weight difference between the current schedule and the MaxWeight schedule, $W^*$ (which is the ``best'' schedule under this measure), and compares this weight difference to a threshold which is a function of the maximum weight, $g(W^*)$. When the schedule weight difference exceeds the threshold, we reconfigure the schedule to the MaxWeight schedule, otherwise keep the current schedule.

The selection of the threshold determines the performance of the policy. In~\cite{AMW_Infocom}, it has been shown that if $g(x) = (1-\gamma)x^{1-\delta}$, then the Adaptive MaxWeight is throughput optimal. The result has been generalized in~\cite{AMW_arxiv} to any sublinear function $g$, and named $g(\cdot)$ as the hysteresis function. The result is stated as the following fact.

\begin{fact}
Given any reconfiguration delay $\Delta_r>0$, and given any sublinear hysteresis function $g$ that is an increasing function, the Markov chain $\mathbf{X}(t)$ is positive recurrent for any admissible traffic rate matrix under the Adaptive MaxWeight. Therefore, the Adaptive MaxWeight is throughput optimal.
\end{fact}

While the throughput optimality is a desirable property, it may be considered as only a first-order performance metric, in the sense that it only guarantees bounded expected queue length (and thus bounded expected delay), but the queue length could still be very large. One more step forward is to characterize its expected queue length. 


\begin{algorithm}
\caption{Adaptive MaxWeight Scheduling Policy}
\begin{algorithmic}
\REQUIRE Sublinear and increasing function $g(\cdot)$
\FOR{each $t = 0, 1, \dots$}
    \STATE $\mathbf{s}^*(t) \leftarrow \arg\max\limits_{\mathbf{s} \in \mathcal{S}} \sum\limits_{ij} q_{ij}(t)s_{ij}(t)$
    \STATE $W^*(t) \leftarrow \max\limits_{\mathbf{s} \in \mathcal{S}} \sum\limits_{ij} q_{ij}(t)s_{ij}(t)$
    \STATE $W(t) \leftarrow \sum\limits_{ij} q_{ij}(t)s_{ij}(t)$
    \STATE $\Delta W(t) \leftarrow W^*(t) - W(t)$ 

    \IF{$\Delta W(t) > g(W^*(t))$}
        \STATE $\mathbf{s}(t+1) \leftarrow \mathbf{s}^*(t)$
    \ELSE
        \STATE $\mathbf{s}(t+1) \leftarrow \mathbf{s}(t)$
    \ENDIF
\ENDFOR
\end{algorithmic}
\label{Alg}

\end{algorithm}

\section{Heavy Traffic Analysis} \label{Sec:heavy}

Studying queue length or delay performance for a queueing system such as a switch in general is challenging. Therefore, such systems are studied in an asymptotic. In this paper, we focus on the heavy traffic regime, and make use of the drift technique developed in~\cite{Atilla_drift}. The outline of the heavy traffic analysis is as follows. We first establish the multi-dimensional state space collapse for the Adaptive MaxWeight. With the state space collapse result, we apply the drift technique to a Lyapunov function proposed in~\cite{MW_scaling}, and obtain an steady state queue length upper bound that is dependent on the expected schedule duration. We then characterize the relation between the expected schedule duration and the queue length, and use this relation to derive bounds on asymptotically tight steady state queue length bounds.

In the following discussion, we consider a sequence of systems indexed by $\epsilon$. Each system has i.i.d. arrival process $\mathbf{a}^{(\epsilon)}(t)$ with mean and variance given by 
\begin{align*}
\bbbe [ \mathbf{a}^{(\epsilon)}(t) ] =& \boldsymbol{\lambda}^{(\epsilon)} = \boldsymbol{\nu} (1 - \epsilon)  , \ \ \ 
\Var [ \mathbf{a}^{(\epsilon)}(t) ] = \big( \boldsymbol{\sigma}^{(\epsilon)} \big)^2
\end{align*}
where $\boldsymbol{\nu} \in \mathcal{F}$ and $\big( \boldsymbol{\sigma}^{(\epsilon)} \big)^2 \rightarrow \boldsymbol{\sigma}^2$ as $\epsilon \rightarrow 0$. The traffic load of each switch is $\rho = 1-\epsilon$. Recall that $\mathcal{F}$ is the set of critically loaded rate matrix with all ports saturated. The sequence of switches considered here have arrival rate matrices that approach  $\boldsymbol{\nu}$ as we take $\epsilon \rightarrow 0$.

\subsection{State Space Collapse}

It was shown in~\cite{MW_scaling} that for a switch with no reconfiguration delay, the MaxWeight scheduling exhibits a multi-dimensional state space collapse. To be specific, let $\mathbf{e}^{(i)}$ denote the matrix with $i^{th}$ row being all ones and zeros everywhere else, and $\tilde{\mathbf{e}}^{(j}$ denote the matrix with $j^{th}$ row being all ones and zeros everywhere else. As $\epsilon \rightarrow 0$, the steady state queue length $\bar{\mathbf{q}}^{(\epsilon)}$ ``concentrates'' in the cone spanned by the matrices $\{\mathbf{e}^{(i)}\}_{i=1}^n \cup \{\tilde{\mathbf{e}}^{(j}\}_{j=1}^n$, \ie 
\begin{align*}
\mathcal{K} = \Big\{  \mathbf{x} \in \bbbr^{n^2}:& \  \mathbf{x} = \sum_i w_i \mathbf{e}^{(i)} + \sum_j \tilde{w}_j \tilde{\mathbf{e}}^{(j)},  
\mbox{ where } w_i, \tilde{w}_j \in \bbbr_+ \mbox{ for all }i, j \Big\},
\end{align*}
in the sense that the projection of $\bar{\mathbf{q}}^{(\epsilon)}$ onto $\mathcal{K}$ is the domniant component in $\bar{\mathbf{q}}^{(\epsilon)}$. More specifically, for any $\mathbf{x} \in \bbbr^{n^2}$, define the projection of $\mathbf{x}$ on to $\mathcal{K}$ as
\[
\mathbf{x}_{\parallel} = \arg\min\limits_{\mathbf{y} \in \mathcal{K}} \|\mathbf{x} - \mathbf{y}\|_2
\]
and define $\mathbf{x}_{\perp} = \mathbf{x} - \mathbf{x}_{\parallel}$. In the heavy traffic limit ($\epsilon \rightarrow 0$), all moments of  $\bar{\mathbf{q}}_{\perp}^{(\epsilon)}$ are bounded by a constant, and hence is a negligible component in $\bar{\mathbf{q}}$ since it can be shown that $\|\bar{\mathbf{q}}\|_1$ is $\Omega(1/\epsilon)$. This is referred as state space collapse (SSC) in~\cite{MW_scaling}.

In this paper, we derive a weaker notion of the state space collapse for switches with reconfiguration delay operated under the Adaptive MaxWeight policy. The following lemma is a $T$-step version of~\cite[Theorem 1]{Bertsimas} where $T$ could be any fixed integer. In the rest of this section, we would use this lemma to derive the (weaker version of) SSC for switches with no reconfiguration.

\begin{lemma}
\label{lemma:Foster}
Consider an irrdeucible and aperiodic Markov Chain $\left\{ X(t) \right\}_{t \geq 0}$ over a countable state space $\mathcal{X}$, suppose $Z: \mathcal{X} \rightarrow \bbbr_+$ is a nonnegative Lyapunov function. For any fixed integer $T > 0$, we define the $T$-step drift $\Delta^T Z(X)$ of $Z$ at state $X$ as
\[
\Delta^T Z(X) = \left[ Z(X(t+T)) - Z(X(t)) \right] \mathds{1}_{\{X(t) = X\}}.
\]
Suppose the $T$-step drift satisfies the following conditions C.1 and C.2:
\begin{itemize}
    \item[C.1] There exists an $\eta > 0$, and a $\kappa < \infty$ such that for any $t = 1,2,\dots$ and for all $X \in \mathcal{X}$ with $Z(X) \geq \kappa$, 
    \[
    \bbbe [ \Delta^T Z(X) | X(t) = X ] \leq -\eta
    \]
    \item[C.2] There exists a $D < \infty$ such that for all $X \in \mathcal{X}$,
    \[
    \Pr \left\{ |\Delta^T Z(X)| \leq D \right\} = 1
    \]
\end{itemize}
If the Markov chain $\left\{ X(t) \right\}_{t \geq 0}$ converges in distribution to a random variable $\bar{X}$, then 
\[
\bbbe [ Z(\bar{X}) ] \leq \kappa + \frac{2D^2}{\eta}.
\]
\end{lemma}

The proof of this lemma could be done by simply replacing the transition probability to $T$-step transition probability in~\cite[Theorem 1]{Bertsimas}, hence we omit the proof here. Lemma~\ref{lemma:Foster} could then be used to proof the following theorem, which is essential to establish the (weak version of) SSC for the Adaptive MaxWeight.

\begin{theorem}
Consider a set of switch systems with a fixed reconfiguration delay $\Delta_r > 0$, parametrized by $0 < \epsilon < 1$, all operated under Adaptive MaxWeight policy with hysteresis function $g(\cdot)$, where $g(\cdot)$ is sublinear and concave. Each system has arrival process $\mathbf{a}^{(\epsilon)}(t)$ as described in Section~\ref{Sec:model}. The mean arrival rate vector $\boldsymbol{\lambda}^{(\epsilon)} = (1-\epsilon) \boldsymbol{\nu}$ for some fixed $\boldsymbol{\nu} \in \mathcal{F}$ such that $\nu_{\min} \stackrel{\Delta}{=} \min\limits_{ij} \nu_{ij} > 0$. Let the variance $\big( \boldsymbol{\sigma}^{(\epsilon)} \big)^2$ of the arrival process satisfy that $\| \boldsymbol{\sigma}^{(\epsilon)} \|^2 \leq \tilde{\sigma}^2$ for some $\tilde{\sigma}^2$ not dependent on $\epsilon$.

Let $\mathbf{X}^{(\epsilon)}(t) \in \mathcal{X}$ denote the process that determines each system, which is positive recurrent and hence converges to a steady state random vector in distribution, denoted as $\bar{\mathbf{X}}^{(\epsilon)} = (\bar{\mathbf{q}}^{(\epsilon)}, \bar{\mathbf{s}}^{(\epsilon)}, \bar{r}^{(\epsilon)})$. Then for any fixed $\theta$ with $0 < \theta < 1/2$, and for each system with $0 < \epsilon \leq \nu_{\min}/4\|\boldsymbol{\nu}\|$, the steady state queue lengths vector satisfies
\[
\bbbe \Big[ \|\bar{\mathbf{q}}^{(\epsilon)}_{\perp}\| - \theta \|\bar{\mathbf{q}}^{(\epsilon)}_{\parallel}\| \Big] \leq M_{\theta}, 
\]
where $M_{\theta}$ is a function of $\theta, \tilde{\sigma}, a_{\max}$, $\nu_{\min}$ and $n$, but is independent of $\epsilon$.
\label{thm:SSC}
\end{theorem}

The proof of Theorem~\ref{thm:SSC} is given in appendix. Comparing Theorem~\ref{thm:SSC} with~\cite[Proposition 2]{MW_scaling}, we may see that we no longer have the guarantee that all moments of $\|\bar{\mathbf{q}}^{(\epsilon)}_{\perp}\|$ are bounded here. However, we could still show that $\bbbe [ \|\bar{\mathbf{q}}^{(\epsilon)}_{\perp}\| ]$ is negligible comparing to $\bbbe [ \|\bar{\mathbf{q}}^{(\epsilon)}\| ]$ as $\epsilon \rightarrow 0$, hence we consider this as a weak version of SSC. In particular, notice that the constant $M_{\theta}$ is independent of $\epsilon$, and that $\bbbe [ \|\bar{\mathbf{q}}^{(\epsilon)}\| ] \rightarrow \infty$ as $\epsilon \rightarrow 0$. Then since $\bbbe [ \|\bar{\mathbf{q}}^{(\epsilon)}_{\perp}\| ] \leq \theta \bbbe [ \|\bar{\mathbf{q}}^{(\epsilon)}_{\parallel}\| ] + M_{\theta} \leq \theta \bbbe [ \|\bar{\mathbf{q}}^{(\epsilon)}\| ] + M_{\theta}$ for any $\epsilon > 0$, we have $\lim\limits_{\epsilon \rightarrow 0} \frac{\bbbe \big[ \|\bar{\mathbf{q}}^{(\epsilon)}_{\perp}\| \big]}{\bbbe \big[ \|\bar{\mathbf{q}}^{(\epsilon)}\| \big] } \leq \theta$ for any $\theta > 0$. Therefore, we may conclude that 
\begin{align}
\lim\limits_{\epsilon \rightarrow 0} \frac{\bbbe \big[ \|\bar{\mathbf{q}}^{(\epsilon)}_{\perp}\| \big]}{\bbbe \big[ \|\bar{\mathbf{q}}^{(\epsilon)}\| \big] }  = 0.
\label{weak_SSC}
\end{align}

\subsection{Drift Analysis}

With the weak SSC result from the previous subsection, we now utilize Lyapunov drift analysis similar to~\cite{MW_scaling} to derive bounds on the steady state of the Markov chain $\mathbf{X}^{(\epsilon)}(t)$. For simplicity of the notation, we drop the superscript $(\epsilon)$ in the following. Now consider the following Lyapunov functions from~\cite{MW_scaling}:
\begin{align*}
V_1(\mathbf{X}) &= \sum_{i} \Big( \sum_{j} q_{ij} \Big)^2, \ \
V_2(\mathbf{X}) = \sum_{i} \Big( \sum_{i} q_{ij} \Big)^2, \ \
V_3(\mathbf{X}) = \Big( \sum_{ij} q_{ij} \Big)^2, \ \ V_4(\mathbf{X}) = V_1(\mathbf{X}) + V_2(\mathbf{X}) - \frac{1}{n} V_3(\mathbf{X})
\end{align*}





It may be shown using Lemma~\ref{lemma:Foster} that for steady state $\bar{\mathbf{X}}$, the expectations $\{ \bbbe[V_k(\bar{\mathbf{X}})] \}_{k=1}^4$ are finite. We thus have zero drift for $V_4$ at steady state:
\begin{align*}
\bbbe_{\bar{\mathbf{X}}} \Big[ \Delta V_4(\mathbf{X}(t)) \Big] 
= \bbbe_{\bar{\mathbf{X}}} \Big[ \Delta V_1(\mathbf{X}(t)) + \Delta V_2(\mathbf{X}(t)) - \frac{1}{n} \Delta V_3(\mathbf{X}(t)) \Big] = 0
\end{align*}

We now evaluate the above drift terms with the queue dynamics and rewrite the expression as
\begin{align*}
\mathcal{T}_1 = \mathcal{T}_2 + \mathcal{T}_3 + \mathcal{T}_4
\end{align*}
where
\begin{align*}
\mathcal{T}_1 
=& \bbbe_{\bar{\mathbf{X}}} \Bigg[  2 \sum_{i} \big(\sum_{j} q_{ij}(t)\big)  \big(\sum_{j}s_{ij}(t)\mathds{1}_{\{r(t)=0\}}- a_{ij}(t) \big) 
+ 2 \sum_{j} \Big(\sum_{i} q_{ij}(t) \Big)  \Big(\sum_{i}s_{ij}(t)\mathds{1}_{\{r(t)=0\}}- a_{ij}(t) \Big) \nonumber \\
&- \frac{2}{n}  \Big(\sum_{ij} q_{ij}(t) \Big)  \Big(\sum_{ij} s_{ij}(t)\mathds{1}_{\{r(t)=0\}}- a_{ij}(t) \Big) \Bigg] \nonumber \\
\mathcal{T}_2
=& \bbbe_{\bar{\mathbf{X}}} \Bigg[ \sum_{i} \Big( \sum_{j} a_{ij}(t) - s_{ij}(t)\mathds{1}_{\left\{r(t)=0\right\}} \Big)^2 
+ \sum_{j} \Big( \sum_{i}a_{ij}(t) - s_{ij}(t)\mathds{1}_{\left\{r(t)=0\right\}} \Big)^2 
- \frac{1}{n} \Big( \sum_{ij}a_{ij}(t) - s_{ij}(t)\mathds{1}_{\left\{r(t)=0\right\}} \Big)^2  \Bigg]  \\
\end{align*}
\begin{align*}
\mathcal{T}_3
=& \bbbe_{\bar{\mathbf{X}}} \Bigg[ - \sum_{i} \Big( \sum_{j}u_{ij}(t) \Big)^2 - \sum_{j} \Big( \sum_{i}u_{ij}(t) \Big)^2  
+ \frac{1}{n} \Big( \sum_{ij}u_{ij}(t) \Big)^2  \Bigg] \\
\mathcal{T}_4
=& \bbbe_{\bar{\mathbf{X}}} \Bigg[ 2 \sum_{i} \Big( \sum_{j}q_{ij}(t+1) \Big) \Big( \sum_{j}u_{ij}(t) \Big)  
+ 2 \sum_{j} \Big( \sum_{i} q_{ij}(t+1) \Big) \Big( \sum_{i}u_{ij}(t) \Big) 
- \frac{2}{n} \Big( \sum_{ij} q_{ij}(t+1) \Big) \Big( \sum_{ij}u_{ij}(t) \Big) \Bigg]
\end{align*}

We now simplify each term:
\begin{align}
\mathcal{T}_1
=& \bbbe_{\bar{\mathbf{X}}} \Big[ 2\Big(\sum_{ij} q_{ij}(t) \Big)  \Big(\mathds{1}_{\{r(t)=0\}}-(1-\epsilon) \Big) \Big]
= \bbbe_{\bar{\mathbf{X}}} \Big[ 2\Big(\sum_{ij} q_{ij}(t) \Big)  \Big(\epsilon - \mathds{1}_{\{r(t)>0\}} \Big)  \Big]  \nonumber \\
=& 2\epsilon \bbbe_{\bar{\mathbf{X}}} \Big[ \sum_{ij} q_{ij}(t)  \Big]  - 2\bbbe_{\bar{\mathbf{X}}} \Big[ \sum_{ij} q_{ij}(t) \mathds{1}_{\{r(t)>0\}} \Big] \nonumber \\
=& 2\epsilon \bbbe_{\bar{\mathbf{X}}} \Big[ \sum_{ij} q_{ij}(t)  \Big]  - 2 \Big( \bbbe_{\bar{\mathbf{X}}} \Big[ \sum_{ij} q_{ij}(t) \Big] - \bbbe_{\bar{\mathbf{X}}} \Big[ \sum_{ij} q_{ij}(t) \Big| r(t)=0 \Big] \Pr{}_{\bar{\mathbf{X}}} \{r(t)=0\} \Big) \nonumber \\
=& 2 \Big( \epsilon - \Pr{}_{\bar{\mathbf{X}}}\{r(t)>0\} \Big) \bbbe_{\bar{\mathbf{X}}} \Big[ \sum_{ij} q_{ij}(t)  \Big]  -  2\Big(\bbbe_{\bar{\mathbf{X}}} \Big[ \sum_{ij} q_{ij}(t) \Big] - \bbbe_{\bar{\mathbf{X}}} \Big[ \sum_{ij} q_{ij}(t) \Big| r(t)=0 \Big] \Big) \Big(1- \Pr{}_{\bar{\mathbf{X}}}\{r(t)>0\} \Big) 
\label{T1}
\end{align}

Note that with the ergodicity of the Markov chain $\mathbf{X}(t)$, $\bbbe_{\bar{\mathbf{X}}} [ \sum_{ij} q_{ij}(t) ]$ and $\bbbe_{\bar{\mathbf{X}}} [ \sum_{ij} q_{ij}(t) | r(t)=0 ]$ equal to the time average of sum queue length and that under $r(t)=0$, respectively. Also note that the probability of $r(t)=0$ approaches $1$ as $\epsilon \rightarrow 0$, which means that the time average under $r(t)=0$ only excludes a diminishing number of time instances. Then since the change in sum queue length $| \sum_{ij}q_{ij}(t+1) - \sum_{ij}q_{ij}(t) | < n^2a_{\max}$ is bounded, we have the difference $\bbbe_{\bar{\mathbf{X}}} [ \sum_{ij} q_{ij}(t) ] - \bbbe_{\bar{\mathbf{X}}} [ \sum_{ij} q_{ij}(t) | r(t)=0 ] \rightarrow 0$ as $\epsilon \rightarrow 0$.

For term $\mathcal{T}_2$, since $\bbbe \Big[ \sum_{i} ( \sum_{j} a_{ij}(t) - 1 )^2 \Big] = \|\boldsymbol{\sigma}\|^2 + n\epsilon^2$ and $\bbbe \Big[ ( \sum_{ij} a_{ij}(t) - n )^2 \Big] = \|\boldsymbol{\sigma}\|^2 + n^2\epsilon^2$, we have

\begin{align}
\mathcal{T}_2 
=& \bbbe_{\bar{\mathbf{X}}} \Bigg[ \Bigg(  \sum_{i} ( \sum_{j} a_{ij}(t) - 1 )^2 + \sum_{j} ( \sum_{i} a_{ij}(t) - 1 )^2 
- \frac{1}{n} ( \sum_{ij} a_{ij}(t) - n )^2 \Bigg)  \nonumber \\
&+ \Bigg(  \sum_{i} ( 2\sum_{j} a_{ij}(t) - 1  ) + \sum_{j} ( 2\sum_{i} a_{ij}(t) - 1  ) 
- \frac{1}{n} ( 2n\sum_{ij} a_{ij}(t) - n^2 ) \Bigg) \mathds{1}_{\{r(t)>0\}} \Bigg]  \nonumber \\
=& \bbbe_{\bar{\mathbf{X}}} \Big[ \Big(  
2 (\|\boldsymbol{\sigma}\|^2 + n\epsilon^2) - \frac{1}{n} (\|\boldsymbol{\sigma}\|^2 + n^2\epsilon^2) \Big) 
+  \Big( 2n(1-\epsilon) - n \Big) \mathds{1}_{\{r(t)>0\}} \Big]  \nonumber \\
=& \Big( (2-\frac{1}{n}) \|\boldsymbol{\sigma}\|^2 + n\epsilon^2 \Big)  +  n(1-2\epsilon) \Pr{}_{\bar{\mathbf{X}}} \{ r(t)>0 \}
\label{T2}
\end{align}

For term $\mathcal{T}_3$, since $u_{ij}(t) \leq s_{ij}(t)$, we have $\sum_{i} u_{ij} \leq 1, \sum_{j} u_{ij} \leq 1$ and $\sum_{ij} u_{ij} \leq n$. Therefore,
\begin{align*}
\mathcal{T}_3 \leq \bbbe_{\bar{\mathbf{X}}} \Bigg[  \frac{1}{n}  \Big( \sum_{ij}u_{ij}(t) \Big)^2 \Bigg] \leq \bbbe_{\bar{\mathbf{X}}} \Bigg[  \sum_{ij}u_{ij}(t)  \Bigg] 
\end{align*}

The above expression is the expected sum of unused services between schedule reconfiguration time instance. One way to determine this value is to set the drift of $\sum_{ij} \bar{\mathbf{q}}_{ij}$ to zero. We may then obtain 
\begin{align}
\bbbe_{\bar{\mathbf{X}}} \Big[ \sum_{ij}u_{ij}(t) \Big] = n \Big( \epsilon - \Pr{}_{\bar{\mathbf{X}}} \{ r(t)>0 \} \Big)
\label{unused_1}
\end{align} 
Hence we have
\begin{align}
\mathcal{T}_3 \leq   n \Big( \epsilon - \Pr{}_{\bar{\mathbf{X}}} \{ r(t)>0 \} \Big)
\label{T3}
\end{align}

We now utilize the relation $u_{ij}(t) q_{ij}(t+1) = 0$ for any $i,j$ to bound the term $\mathcal{T}_4$. The expression below follows along the same lines in~\cite{MW_scaling}:
\begin{align*}
\mathcal{T}_4 =& 2 \bbbe_{\bar{\mathbf{X}}} \Bigg[  \sum_{ij} u_{ij}(t) \Big( \sum_{j'} q_{ij'}(t+1) + \sum_{i'} q_{i'j}(t+1) 
- \frac{1}{n} \sum_{i'j'} q_{i'j'}(t+1) \Big) \Bigg] \\
=& 2 \bbbe_{\bar{\mathbf{X}}} \Big[  \Big\langle \mathbf{u}(t), -n \mathbf{q}_{\perp}(t+1) + \sum_{i} \langle \mathbf{q}_{\perp}(t+1), \mathbf{e}^{(i)} \rangle \mathbf{e}^{(i)} 
+ \sum_{j} \langle \mathbf{q}_{\perp}(t+1), \tilde{\mathbf{e}}^{(j)} \rangle \tilde{\mathbf{e}}^{(j)}  - \frac{1}{n} \langle \mathbf{q}_{\perp}(t+1), \mathbf{1} \rangle \mathbf{1} \Big\rangle \Big] \\
\leq& 2 \bbbe_{\bar{\mathbf{X}}} \Big[  \Big\langle \mathbf{u}(t), -n \mathbf{q}_{\perp}(t+1) - \frac{1}{n} \langle \mathbf{q}_{\perp}(t+1), \mathbf{1} \rangle \mathbf{1} \Big\rangle \Big]  \\
\leq& 2 \bbbe_{\bar{\mathbf{X}}} \Big[  \big\| \mathbf{u}(t) \big\|  \big\| -n \mathbf{q}_{\perp}(t+1) - \frac{1}{n} \langle \mathbf{q}_{\perp}(t+1), \mathbf{1} \rangle \mathbf{1} \big\|  \Big] 
\end{align*}
where the last inequality follows from Cauchy-Schwartz inequality.

Note that 
\begin{align*}
\big\| -n \mathbf{q}_{\perp}(t+1) - \frac{1}{n} \langle \mathbf{q}_{\perp}(t+1), \mathbf{1} \rangle \mathbf{1} \big\| 
\stackrel{(a)}{\leq}& n \| \mathbf{q}_{\perp}(t+1) \| + \frac{1}{n} |\langle \mathbf{q}_{\perp}(t+1), \mathbf{1} \rangle| \| \mathbf{1} \| \\
\stackrel{(b)}{\leq}& n \| \mathbf{q}_{\perp}(t+1) \| + \frac{\| \mathbf{1} \| \| \mathbf{1} \| }{n} \| \mathbf{q}_{\perp}(t+1)\| \\
=& 2n \| \mathbf{q}_{\perp}(t+1) \|
\end{align*}
where $(a)$ follows from triangle inequality, and $(b)$ follows from Cauchy-Schwartz inequality.

We then obtain
\begin{align*}
\mathcal{T}_4 
\leq& 4n \bbbe_{\bar{\mathbf{X}}} \Big[  \|\mathbf{u}(t)\| \|\mathbf{q}_{\perp}(t+1) \| \Big]  \\
\leq& 4n \bbbe_{\bar{\mathbf{X}}} \Big[  \Big( \sum_{ij} u_{ij}(t) \Big) \Big( \|\mathbf{q}_{\perp}(t)\| + na_{\max} \Big) \Big]  \\
=& 4n \bbbe_{\bar{\mathbf{X}}} \Big[ \Big( \sum_{ij} u_{ij}(t) \Big) \|\mathbf{q}_{\perp}(t)\| \Big] +  4n \bbbe_{\bar{\mathbf{X}}} \Big[  \sum_{ij} u_{ij}(t) \Big] na_{\max} 
\end{align*}

Since $\sum_{ij}u_{ij}(t) \leq n$, we may write $\sum_{ij}u_{ij}(t) \leq n\mathds{1}_{\{\sum_{ij}u_{ij}(t) > 0\}}$ and thus
\begin{align*}
&\bbbe_{\bar{\mathbf{X}}} \Big[ \Big( \sum_{ij} u_{ij}(t) \Big) \|\mathbf{q}_{\perp}(t)\| \Big] \\
\leq& n \bbbe_{\bar{\mathbf{X}}} \Big[ \mathds{1}_{\{\sum_{ij} u_{ij}(t) > 0 \}} \|\mathbf{q}_{\perp}(t)\| \Big]  \\
=& n \Big( \bbbe_{\bar{\mathbf{X}}} \Big[ \|\mathbf{q}_{\perp}(t)\| \Big] - 
\bbbe_{\bar{\mathbf{X}}} \Big[ \|\mathbf{q}_{\perp}(t)\| \mathds{1}_{\{\sum_{ij} u_{ij}(t) = 0 \}} \Big]  \Big) \\
=& n \Big( \bbbe_{\bar{\mathbf{X}}} \Big[ \|\mathbf{q}_{\perp}(t)\| \Big] - 
\bbbe_{\bar{\mathbf{X}}} \Big[ \|\mathbf{q}_{\perp}(t)\|  \Big| \sum_{ij} u_{ij}(t) = 0  \Big] \Pr{}_{\bar{\mathbf{X}}} \Big\{\sum_{ij} u_{ij}(t)=0\Big\}  \Big) \\
=& n \bbbe_{\bar{\mathbf{X}}} \Big[ \|\mathbf{q}_{\perp}(t)\| \Big] \Pr{}_{\bar{\mathbf{X}}} \Big\{\sum_{ij} u_{ij}(t)>0\Big\} + 
n \Big( \bbbe_{\bar{\mathbf{X}}} \Big[ \|\mathbf{q}_{\perp}(t)\| \Big] - \bbbe_{\bar{\mathbf{X}}} \Big[ \|\mathbf{q}_{\perp}(t)\|  \Big| \sum_{ij} u_{ij}(t) = 0  \Big] \Big) \Pr{}_{\bar{\mathbf{X}}} \Big\{\sum_{ij} u_{ij}(t)=0\Big\} 
\end{align*}

From (\ref{unused_1}), we have $\bbbe_{\bar{\mathbf{X}}} [  \sum_{ij} u_{ij}(t) ] = n(\epsilon - \Pr{}_{\bar{\mathbf{X}}} \{r(t)>0\} ) = \sum_{k=1}^{n}  \Pr_{\bar{\mathbf{X}}} \{ \sum_{ij} u_{ij}(t) = k\} \cdot k \geq  \sum_{k=1}^{n} \Pr_{\bar{\mathbf{X}}} \{ \sum_{ij} u_{ij}(t) = k\} = \Pr_{\bar{\mathbf{X}}} \{ \sum_{ij} u_{ij}(t) > 0\}$, hence 
\begin{align}
\mathcal{T}_4 
\leq 4n^{3} \Big( \epsilon - \Pr{}_{\bar{\mathbf{X}}} \{r(t)>0\} \Big) \Big( \bbbe_{\bar{\mathbf{X}}} \Big[ \|\mathbf{q}_{\perp}(t)\| \Big] + a_{\max} \Big) + 4n^{2}\Big( \bbbe_{\bar{\mathbf{X}}} \Big[ \|\mathbf{q}_{\perp}(t)\| \Big] - \bbbe_{\bar{\mathbf{X}}} \Big[ \|\mathbf{q}_{\perp}(t)\|  \Big| \sum_{ij} u_{ij}(t) = 0  \Big] \Big) 
\label{T4}
\end{align}

Similar to the second term in $\mathcal{T}_1$, by the ergodicity of the Markov chain $\mathbf{X}(t)$, and that the probability of $\sum_{ij}u_{ij}(t)=0$ approaches $1$ as $\epsilon \rightarrow 0$. Then since the change in $\|\mathbf{q}_{\perp}(t)\|$, \textsl{i.e.} $\big| \|\mathbf{q}_{\perp}(t+1)\| - \|\mathbf{q}_{\perp}(t)\| \big| < na_{\max}$ is bounded, we have that $\bbbe_{\bar{\mathbf{X}}} [ \|\mathbf{q}_{\perp}(t)\| ] - \bbbe_{\bar{\mathbf{X}}} [ \|\mathbf{q}_{\perp}(t)\| | \sum_{ij}u_{ij}(t)=0 ] \rightarrow 0$ as $\epsilon \rightarrow 0$.

Combining (\ref{T1})-(\ref{T4}), we obtain
\begin{align}
\Big( \epsilon - \Pr{}_{\bar{\mathbf{X}}}\{r(t)>0\} \Big) \Big( \bbbe_{\bar{\mathbf{X}}} \Big[ \sum_{ij} q_{ij}(t)  \Big] - 2n^{3} \bbbe_{\bar{\mathbf{X}}} \Big[ \|\mathbf{q}_{\perp}(t)\| \Big] \Big)
\leq 
(1-\frac{1}{2n}) \|\boldsymbol{\sigma}\|^2 + B(\epsilon, n)
\label{queue_upper_bound}
\end{align}
where $\lim_{\epsilon \rightarrow 0} B(\epsilon, n) = 0$.


\subsection{Expected Schedule Duration}


In the previous subsection, we derived a bound on steady state queue lengths. Note that this bound depends on the probability of the switch in reconfiguration $\Pr{}_{\bar{\mathbf{X}}}\{r(t) > 0\}$. In this subsection, we derive the mean schedule duration in order to evaluate this probability.

Recall that in the previous subsection, we set the drift of $\sum_{ij} \bar{q}_{ij}$ to zero to obtain (\ref{unused_1}), which is an expression of the total unused service. Here we consider the drift of another Lyapunov function to obtain a different expression for the total unused service, and combined the two expressions to derive the expected schedule duration. 

In this subsection, we consider the original Markov chain $\mathbf{X}^{(\epsilon)}(t)$ sampled at the reconfiguration times $\{t_k^S\}$ and denote it as $\mathbf{X}^{(\epsilon)}_k = \mathbf{X}^{(\epsilon)}(t_k^S)$. Note that $\{t_k^S\}$ is a stopping time with respect to $\mathbf{X}^{(\epsilon)}(t)$, hence $\mathbf{X}^{(\epsilon)}_k$ is also a Markov chain. 

Define the Lyapunov function $W$ on state $\mathbf{X} = (\mathbf{q}, \mathbf{s}, r)$:
\begin{align*}
W(\mathbf{X}) = \Big\langle \mathbf{q}, \mathbf{s} \Big\rangle = \sum_{ij} q_{ij} s_{ij}
\end{align*}
which is simply the schedule weight function. Note that $W(\mathbf{X}) \leq \sum_{ij} q_{ij}$, hence the steady state mean of $W(\mathbf{X})$ is finite. We may then set the drift of $W(\mathbf{X})$ between two schedule reconfiguration time instance to be zero:
\begin{align*}
&\bbbe \Bigg[ \sum_{t=t_k}^{t_{k+1}-1} \Delta W(\mathbf{X}(t)) \Bigg| \mathbf{X}(t_k) = \hat{\mathbf{X}} \Bigg]  \\
=& \bbbe \Bigg[ \sum_{t=t_k}^{t_{k+1}-1} \Big( \sum_{ij} q_{ij}(t+1)s_{ij}(t+1) - \sum_{ij} q_{ij}(t)s_{ij}(t) \Big)  \Bigg| \mathbf{X}(t_k) = \hat{\mathbf{X}} \Bigg]  \nonumber \\
\stackrel{(a)}{=}& \bbbe \Bigg[ \sum_{ij} q_{ij}(t_{k+1}) \Big( s_{ij}(t_{k+1}) - s_{ij}(t_k) \Big) 
+ \sum_{t=t_k}^{t_{k+1}-1}  \sum_{ij} \Big(q_{ij}(t+1) - q_{ij}(t) \Big) s_{ij}(t_k) \Bigg| \mathbf{X}(t_k) = \hat{\mathbf{X}} \Bigg]  \nonumber \\
\stackrel{(b)}{=}& \bbbe \Bigg[ g(W^*(t_{k+1})) + \delta_{W} + \sum_{t=t_k}^{t_{k+1}-1}  \Big( \sum_{ij} \lambda_{ij}s_{ij}(t_k) - n\mathds{1}_{\{r(t)=0\}} \Big) 
+ \sum_{t=t_k}^{t_{k+1}-1} \sum_{ij} u_{ij}(t)  \Bigg| \mathbf{X}(t_k) = \hat{\mathbf{X}} \Bigg]  \\
=& 0
\end{align*}
Since the schedule remains $\mathbf{s}(t_k)$ between $t_k$ and $t_{k+1}$, and changes to $\mathbf{s}(t_{k+1})$ at time $t_{k+1}$, we have (a). Then by the definition of the Adaptive MaxWeight, the weight difference exceeds $g(W^*(t_{k+1}))$ at time $t_{k+1}$ and thus we may write $\sum_{ij} q_{ij}(t_{k+1}) \Big( s_{ij}(t_{k+1}) - s_{ij}(t_k) \Big) = g(W^*(t_{k+1})) + \delta_{W}$, where $0 < \delta_{W} < na_{\max}$, we have (b).

Note $\bbbe \big[ \sum_{ij} \lambda_{ij}s_{ij}(t_k) \big| \mathbf{X}(t_k) = \hat{\mathbf{X}} \big] = (1-\epsilon) \big\langle \boldsymbol{\nu}, \bbbe[ \mathbf{s}(t_k) | \mathbf{X}(t_k) = \hat{\mathbf{X}} ] \big\rangle$, we may then just denote $\big\langle \boldsymbol{\nu}, \bbbe[ \mathbf{s}(t_k) | \mathbf{X}(t_k) = \hat{\mathbf{X}} ] \big\rangle$ as $\alpha$ and get
\begin{align}
\bbbe \Bigg[ \sum_{t=t_k}^{t_{k+1}-1}  \sum_{ij}  u_{ij}(t) \Bigg| \mathbf{X}(t_k) = \hat{\mathbf{X}} \Bigg]  
=& \Big( n - \alpha (1-\epsilon) \Big)  \bbbe \Big[ t_{k+1}-t_k \Big| \mathbf{X}(t_k) = \hat{\mathbf{X}} \Big]  
- n\Delta_r - \bbbe \Big[ g(\hat{\mathbf{W}}^*) + \delta_W \Big]
\label{unused_2}
\end{align}

Similar to (\ref{unused_1}), we set the drift of $\sum\limits_{ij} q_{ij}$ to zero and get $\bbbe \Big[ \sum\limits_{t=t_k}^{t_{k+1}-1} \sum\limits_{ij} u_{ij}(t) \Big| \mathbf{X}(t_k) = \hat{\mathbf{X}} \Big] = n \Big( \epsilon \bbbe \Big[ t_{k+1} - t_k \Big| \mathbf{X}(t_k) = \hat{\mathbf{X}} \Big] - \Delta_r \Big)$. Combine this with (\ref{unused_2}), we then have
\begin{align}
\bbbe \Big[ t_{k+1}-t_k \Big| \mathbf{X}(t_k) = \hat{\mathbf{X}} \Big] 
= \frac{ \bbbe \Big[ g(\hat{\mathbf{W}}^*) + \delta_W  \Big] }{ (n-\alpha)(1-\epsilon)}  
\label{expected_T}
\end{align}
This essentially connects the expected schedule duration to the expected queue length.

From (\ref{unused_1}), we have 
\begin{align*}
\epsilon \bbbe \Big[ t_{k+1}-t_k \Big| \mathbf{X}(t_k) = \hat{\mathbf{X}} \Big] - \Delta_r = \frac{1}{n} \bbbe_{\hat{\mathbf{X}}} \Bigg[ \sum_{t=t_k}^{t_{k+1}-1}  \sum_{ij}  u_{ij}(t) \Bigg| \mathbf{X}(t_k) = \hat{\mathbf{X}} \Bigg] \geq 0
\end{align*}

Together with (\ref{expected_T}), we have a lower bound on the expected schedule duration, which then implies a lower bound on the expected maximum weight through Jensen's inequality as follows:
\begin{align}
&\bbbe_{\hat{\mathbf{X}}} \Big[ t_{k+1}-t_k \Big| \mathbf{X}(t_k) = \hat{\mathbf{X}} \Big] 
\geq \frac{ \Delta_r}{\epsilon}  \\
\Rightarrow &\bbbe \Big[ \sum_{ij} \hat{\mathbf{q}}_{ij} \Big] \geq
\bbbe \Big[ \hat{\mathbf{W}}^* \Big] 
\geq g^{-1}\Big( \bbbe \Big[ g(\hat{\mathbf{W}}^*) \Big] \Big) 
\geq g^{-1} \bigg( \frac{(n-\alpha)(1-\epsilon) \Delta_r}{\epsilon} - \bbbe \Big[ \delta_W \Big] \bigg)
\label{weight_lower_bound}
\end{align}

Since $\bbbe \Big[ \delta_W \Big] < na_{\max}$ by definition, hence in the heavy traffic regime ($\epsilon \downarrow 0$), we have $\bbbe \Big[ \sum_{ij} \hat{\mathbf{q}}_{ij} \Big] \sim \Omega(g^{-1}(1/\epsilon))$.

By the ergodicity of the Markov chain, we may use the expected schedule duration to derive the probability that the switch is in reconfiguration delay as 
\begin{align}
\Pr{}_{\bar{X}}\{ r(t) > 0 \} = \bbbe_{\bar{\mathbf{X}}} \Big[ \mathds{1}_{\{r(t)>0\}} \Big] = \frac{\Delta_r}{\bbbe [ t_{k+1}-t_k | \mathbf{X}(t_k) = \hat{\mathbf{X}} ] } = \frac{(n-\alpha)(1-\epsilon)\Delta_r}{\bbbe [ g(\hat{\mathbf{W}}^*) + \delta_W ]}
\label{prob_reconfiguration}
\end{align}

We may then apply (\ref{prob_reconfiguration}) into (\ref{queue_upper_bound}). For the simplicity of notation, we denote $\hat{\mathbf{g}} = \bbbe \Big[ g(\hat{W}^*) + \delta_W \Big] $, $\bar{\mathbf{q}}_{s} = \bbbe_{\bar{\mathbf{X}}} \Big[   \sum_{ij} \bar{q}_{ij} \Big] - 4n^3\bbbe_{\bar{\mathbf{X}}} \Big[ \|\bar{q}_{\perp}\| \Big]$, $\beta = (1 - \frac{1}{2n}) \tilde{\sigma}^2 + B(\epsilon, n)$ and obtain

\begin{align*}
&\Big( \epsilon \hat{\mathbf{g}} - (n-\alpha)(1-\epsilon)\Delta_r \Big) \bar{\mathbf{q}}_{s} \leq \beta \hat{\mathbf{g}} \\
\Rightarrow \ & \Big( \epsilon \bar{\mathbf{q}}_{s} - \beta \Big) \Big( \hat{\mathbf{g}} + (n-\alpha) \Delta_r \Big) \leq  (n-\alpha) \Delta_r (\bar{\mathbf{q}}_{s} - \beta)
\end{align*}

We thus have 
\begin{align*}
\hat{\mathbf{g}}
\leq (n-\alpha) \Delta_r \Bigg( \frac{\bar{\mathbf{q}}_{s} - \beta}{\epsilon \bar{\mathbf{q}}_{s} - \beta} - 1 \Bigg) 
= (n-\alpha) \Delta_r \frac{1-\epsilon}{\epsilon}  \Big(  1 + \frac{\beta}{\epsilon \bar{\mathbf{q}}_{s} - \beta } \Big)
\end{align*}

Note that from the lower bound (\ref{weight_lower_bound}) and the weak SSC result~(\ref{weak_SSC}), we have $\bar{\mathbf{q}}_{s} \sim \Omega \Big( g^{-1}(1/\epsilon) \Big)$. Therefore $\frac{\beta}{\epsilon \bar{\mathbf{q}}_{s} - \beta } \rightarrow 0$ as $\epsilon \downarrow 0$, and thus
\begin{align}
&\lim\sup\limits_{\epsilon \downarrow 0} \  \hat{\mathbf{g}}  = \frac{(n-\alpha) \Delta_r}{\epsilon}  
\label{asymptotic_upper_bound}
\end{align}

Along with the lower bound from (\ref{weight_lower_bound}), and that $\bbbe[\delta_W] < na_{\max}$, we have 
\begin{align}
\lim_{\epsilon \downarrow 0} \bbbe \Big[ g(\hat{\mathbf{W}}^*) \Big] = \frac{(n-\alpha)\Delta_r}{\epsilon}
\label{scaling_g}
\end{align}

Note that the scaling in~(\ref{scaling_g}) does not directly imply the scaling of the expected queue length, and Jensen's inequality could only be used to derive the $g^{-1}(1/\epsilon)$ asymptotic lower bound. However, it may be desirable to approximate the scaling of the sum of queue length $\bbbe \Big[ \sum_{ij} \hat{\mathbf{q}}_{ij} \Big]$ with $O(n g^{-1}(1/\epsilon))$, which becomes more accurate as the hysteresis function $g$ becomes closer to linear (\textsl{i.e.} $\delta\rightarrow 0$). In fact, this is also the regime of interest for optimal delay scaling, since $O(n g^{-1}(1/\epsilon))$ scaling improves as $g$ becomes closer to linear. In other words, as we take $\delta \rightarrow 0$ for $g(x) = (1-\gamma)x^{1-\delta}$, we not only get a tighter asymptotic bounds but also a better scaling. One caveat here is that selecting $g$ as exactly a linear function does not fit the analysis in this paper. In fact, it is still even unclear whether throughput optimal could be guaranteed if $g$ is linear.

\section{Simulations} \label{Sec:simulation}

In this section, we show simulation results for switches with reconfiguration delay operated under the Adptive MaxWeight policy, with hysteresis function $g(x) = (1-\gamma) x^{1-\delta}$. We compare the simulation result with the heavy traffic queue length bound obtained in the previous section. 

We now briefly describe the simulation setup. The arrival processes are assumed to be Poisson processes, all with the same arrival rate, which is also known as the uniform traffic. More specifically, the matrix $\boldsymbol{\nu} \in \mathcal{F}$ satisfies $\nu_{ij} = \frac{1}{n}$, $\forall i, j \in \{1, \dots, n\}$. Therefore, we have $\|\boldsymbol{\nu}\|^2 = 1$, and thus $\alpha = n-1$ in  (\ref{asymptotic_upper_bound}). For the parameter of the hysteresis function $g$, since we are only interested in the scaling, we fix $\gamma = 0.1$, and consider $\delta \in \{0.05, 0.1, 0.2\}$ for mean queue length comparison.

\begin{figure}[!t]
\centering
\includegraphics[height=2.1in]{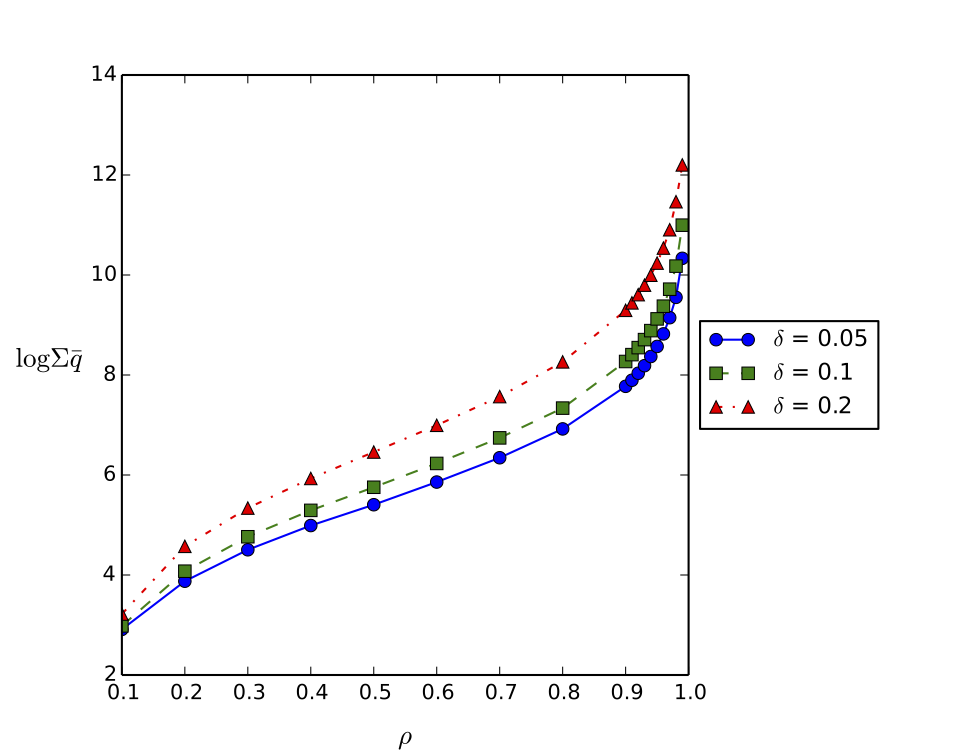}
\caption{Mean total queue length versus traffic load $\rho$ under uniform traffic. Number of ports is $n = 4$, and reconfiguration delay $\Delta_r = 20$.
\vspace{-0.3in}
}
\label{logQ_rho}
\end{figure}

\begin{figure}[!t]
\centering
\includegraphics[height=2.0in]{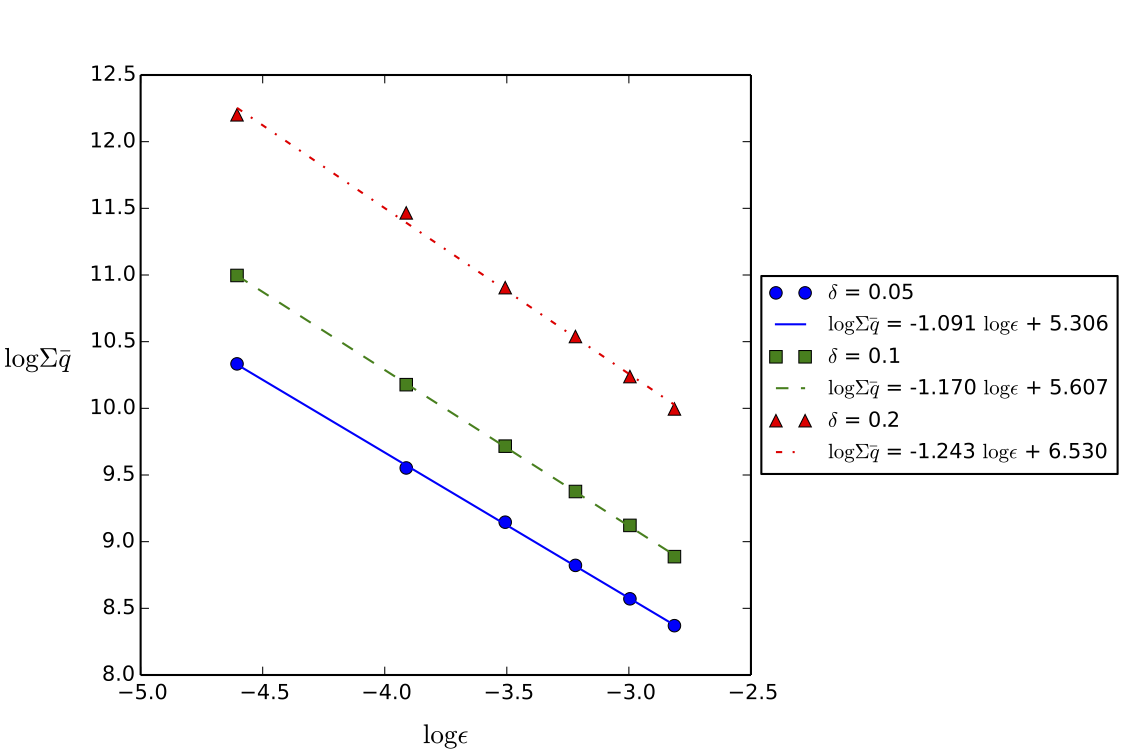}
\caption{Mean total queue length versus traffic load near the capacity region. Number of ports is $n = 4$, and reconfiguration delay $\Delta_r = 20$.
\vspace{-0.3in}
}
\label{logQ_logE}
\end{figure}

Fig.~\ref{logQ_rho} shows the average queue length under traffic loads $\rho \in [0.1, 1]$. The purpose of Fig.~\ref{logQ_rho} is to illustrate that while our analysis focus on critical traffic loads, the conclusion that $\delta$ close to zero gives better scaling also applies for lower traffic loads. In Fig.~\ref{logQ_logE}, we plot the queue length for $\epsilon \in [0.01, 0.06]$ in log scale (corresponding to $\rho \in [0.94, 0.99]$). We may then use linear regression to determine the scaling (\ie the exponent) of mean queue length in the heavy traffic regime. With the scaling result from (\ref{asymptotic_upper_bound}), the scaling with respect to $\epsilon$ is close to $g^{-1}(1/\epsilon)$, hence the theoretical exponent should be $- 1/(1-\delta)$, which would be $\{-1.053, -1.111, -1.250\}$ for $\delta = \{0.05, 0.1, 0.2\}$, respectively.

\begin{figure}[!t]
\centering
\includegraphics[height=2.2in]{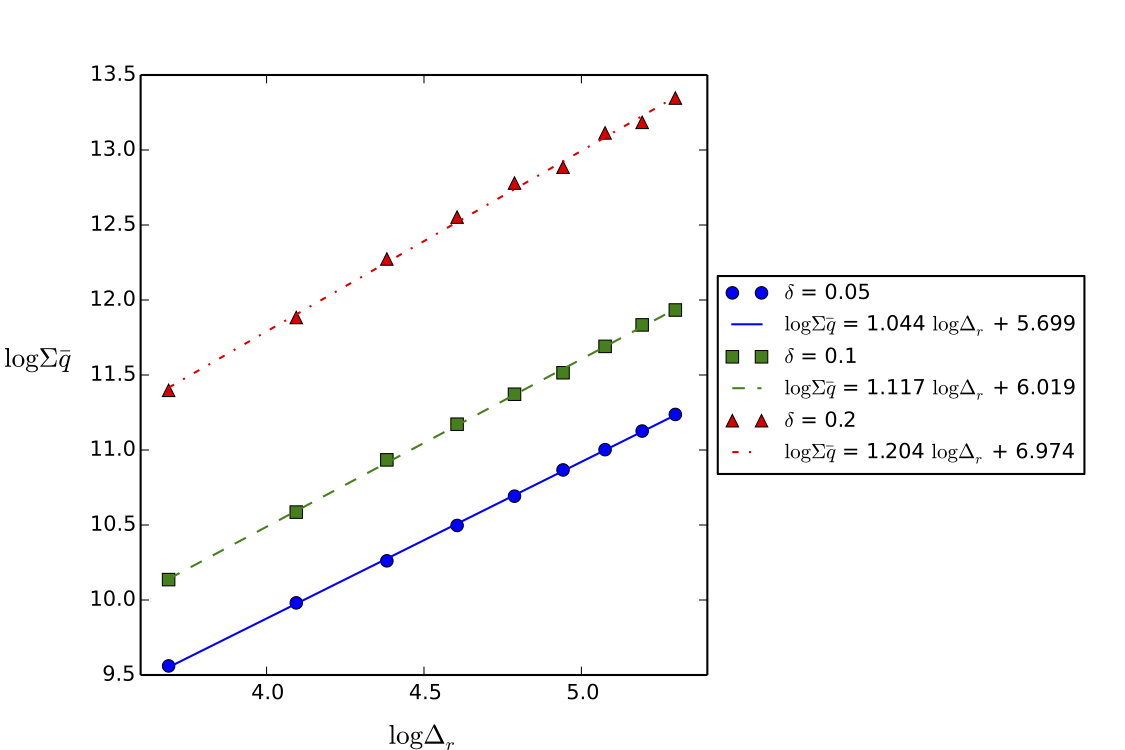}
\caption{Mean total queue length versus reconfiguration delay $\Delta_r$. Number of ports is $n = 4$, and traffic load is $\rho = 0.96$.
\vspace{-0.3in}
}
\label{logQ_logDeltar}
\end{figure}

\begin{figure}[!t]
\centering
\includegraphics[height=2.2in]{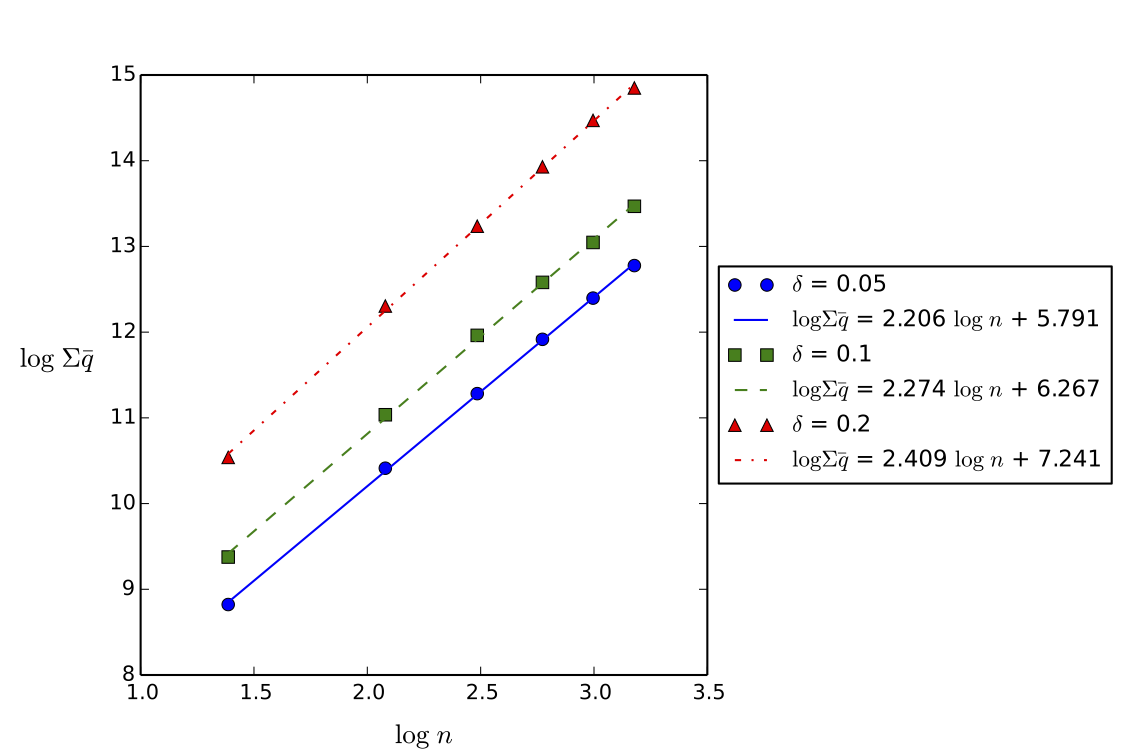}
\caption{Mean total queue length versus number of ports $n$. Reconfiguration delay is $\Delta_r = 20$, and traffic load is $\rho = 0.96$.
\vspace{-0.3in}
}
\label{logQ_logN}
\end{figure}

Figs.~\ref{logQ_logDeltar}~and~\ref{logQ_logN} show the queue length scaling behavior under varying reconfiguration delay $\Delta_r$ and varying number of ports $n$, while the traffic load is fixed as $\rho = 0.96$ (or $\epsilon = 0.04$). For the reconfiguration delay, the scaling is $g^{-1}(\Delta_r)$, hence the theoretical exponents are $\{1.053, 1.111, 1.250\}$ for $\delta \in \{0.05, 0.1, 0.2\}$, respectively. We may see from Fig.~\ref{logQ_logDeltar} that the exponents obtained from the simulation result are very close to our derived scaling. On the other hand, the scaling with respect to $n$ is $n g^{-1}(n)$, hence the theoretical exponents should be $\{2.053, 2.111, 2.250\}$. We could see that the exponents derived from the simulation result are slightly larger than the our derived scaling. One possible reason is that $\epsilon$ is not close enough to $0$, hence there are higher exponent terms that are diminishing in heavy traffic regime affect the scaling.

\section{Conclusion} \label{Sec:conclusion}

We consider the heavy traffic queue length behavior in an input-queued switch with reconfiguration delay, operating under the Adaptive MaxWeight policy. It is shown that the Adaptive MaxWeight exhibits weak state space collapse behavior, which could be considered as an inheritance from the MaxWeight policy in the regime of zero reconfiguration delay. Utilizing the Lyapunov drift technique introduced in~\cite{MW_scaling}, we obtain a queue length upper bound in heavy traffic, which depends on the expected schedule duration. We then discover a relation between the expected schedule duration and the expected queue length, which then implies asymptotically tight bounds for the expected schedule duration in heavy traffic limit, thus determining its scaling. The scaling of the expected schedule duration then implies the dependence of the queue length scaling with the selection of the hysteresis function $g$, and that this scaling improves as $g$ becomes closer to linear. Simulation results are also presented to illustrate the queue length scaling with respect to several system parameters (e.g. traffic load, number of ports, reconfiguration delay) and for comparison to the derived queue length scaling in heavy traffic.

The results obtained in this paper apply to traffic patterns that all input and output ports are saturated. It would be interested to consider the queue length behavior of the Adaptive MaxWeight under incompletely saturated traffic, for example, traffic conditions considered in~\cite{MW_scaling_incompletely_saturated}. The state space collapse result might be similar due to the inheritance from MaxWeight policy, but the characterization for the expected schedule duration remains unclear at this point.

Another interesting aspect for the analysis of queue length behavior of switches with reconfiguration delay is that whether there exists a policy that achieves the optimal scaling. The queue length scaling result presented in this paper suggests that the hysteresis function $g$ should be as close to linear as possible in order to achieve good heavy traffic scaling. However, the throughput analysis currently known excludes the case that $g$ being linear, and thus the analysis in this paper does not apply. This result leaves the existence of an optimal scaling policy as an open problem.

\bibliographystyle{ieeetr}
\bibliography{reference}

\appendix

\begin{proof} (of Theorem~\ref{thm:SSC}) \\
For ease of notation, we drop the superscript $(\epsilon)$ in the following derivation. For each state $\mathbf{X} = (\mathbf{q}, \mathbf{s}, r)$, we define the Lyapunov function $Z(\mathbf{X}) = \max\{ \|\mathbf{q}_{\perp}\| - \theta \| \mathbf{q}_{\parallel} \|, 0 \}$. We then apply Lemma~\ref{lemma:Foster} with the Lyapunov function $Z$ to obtain the result. Note that the selection of the Lyapunov function is such that $Z$ is a nonnegative function. Since $\|\mathbf{q}_{\perp}\| - \theta \| \mathbf{q}_{\parallel} \| \leq Z(\mathbf{X})$ for any state $\mathbf{X} = (\mathbf{q}, \mathbf{s}, r)$, the result follows a bound on $\bbbe[Z(\bar{\mathbf{X}})]$.

We first verify Condition C.2 for $Z(\mathbf{X})$:
\begin{align}
|\Delta^T Z(\mathbf{X})| 
=& \Big| \Big(\|\mathbf{q}_{\perp}(t+T)\| - \theta \|\mathbf{q}_{\parallel}(t+T)\| \Big) - \Big(\|\mathbf{q}_{\perp}(t)\| - \theta \|\mathbf{q}_{\parallel
}(t)\| \Big) \Big| \nonumber \\
\leq& \Big| \|\mathbf{q}_{\perp}(t+T)\| - \|\mathbf{q}_{\perp}(t)\| \Big| + \theta \Big| \|\mathbf{q}_{\parallel}(t+T)\| - \|\mathbf{q}_{\parallel}(t)\| \Big|  \nonumber \\
\leq& \|\mathbf{q}_{\perp}(t+T) - \mathbf{q}_{\perp}(t)\| + \theta  \|\mathbf{q}_{\parallel}(t+T) - \mathbf{q}_{\parallel}(t)\|  \nonumber \\
\leq& (1+\theta) \left\| \mathbf{q}(t+T) - \mathbf{q}(t) \right\|  \nonumber \\
\leq& (1+\theta) na_{\max}T  \triangleq D
\label{drift_deterministic_bound}
\end{align}
Here we use the fact that $\mathbf{q}_{\perp}$ is a projection onto $\mathcal{K}^{\circ} = \{\mathbf{x} \in \bbbr^{n^2}: \langle \mathbf{x}, \mathbf{y} \rangle,  \forall \mathbf{y} \in \mathcal{K} \}$, the polar cone of $\mathcal{K}$. Since the projection onto a cone is non-expansive, we have $\| \mathbf{x}_{\perp} - \mathbf{y}_{\perp} \| \leq \|\mathbf{x} - \mathbf{y}\|$ and $\| \mathbf{x}_{\parallel} - \mathbf{y}_{\parallel} \| \leq \|\mathbf{x} - \mathbf{y}\|$, $\forall \mathbf{x}, \mathbf{y}$.

To verify condition C.1, we need to bound the $T$-step drift for $Z(\mathbf{X})$. For ease of notation, we denote $\bbbe[\ \cdot \ |\mathbf{X}(t) = \mathbf{X}]$ as $\bbbe_{\mathbf{X}}[\ \cdot\ ]$. From~(\ref{drift_deterministic_bound}), it is not hard to see that $\forall \mathbf{X}: Z(\mathbf{X}) > D$,
\begin{align}
\bbbe_{\mathbf{X}}\Big[\Delta^TZ(\mathbf{X}) \Big] = \bbbe_{\mathbf{X}} \Big[ \Big(\|\mathbf{q}_{\perp}(t+T)\|-\|\mathbf{q}_{\perp}(t)\| \Big) - \theta \Big(\|\mathbf{q}_{\parallel}(t+T)\|-\|\mathbf{q}_{\parallel}(t)\| \Big)  \Big].
\label{drift_Z}
\end{align}
Therefore, we need only to consider the $T$-step expected drift of $\|\mathbf{q}_{\perp}\|$ and $\|\mathbf{q}_{\parallel}\|$.

We first consider the drift of $\|\mathbf{q}_{\perp}\|$. The derivation follows the line in \cite{MW_scaling} where the relation in~\cite[Lemma 4]{MW_scaling} is used: Let $V(\mathbf{X}) = ||\mathbf{q}||^2$, $V_{\parallel}(\mathbf{X}) = ||\mathbf{q}_{\parallel}||^2$, and $\Delta V$, $\Delta V_{\parallel}$ denote the one-step drift of $V$, $V_{\parallel}$ (respectively), then $\|\mathbf{q}_{\perp}(t+1)\| - \|\mathbf{q}_{\perp}(t)\|  \leq \frac{1}{2\|\mathbf{q}_{\perp}(t)\|} \big(\Delta V(\mathbf{X}(t)) - \Delta V_{\parallel}(\mathbf{X}(t)) \big)$. With the relation, we have
\begin{align}
\bbbe_{\mathbf{X}} \Big[ \|\mathbf{q}_{\perp}(t+T)\|-\|\mathbf{q}_{\perp}(t)\|  \Big]
=& \bbbe_{\mathbf{X}} \bigg[ \sum_{\tau=t}^{t+T-1} \big( \|\mathbf{q}_{\perp}(\tau+1)\| - \|\mathbf{q}_{\perp}(\tau) \| \big) \bigg]  \nonumber \\
\leq& \bbbe_{\mathbf{X}} \bigg[ \sum_{\tau=t}^{t+T-1} \frac{\Delta V(\mathbf{X}(\tau)) - \Delta V_{\parallel}(\mathbf{X}(\tau))}{2\|\mathbf{q}_{\perp}(\tau)\|} \bigg]  \nonumber \\
=& \bbbe_{\mathbf{X}} \bigg[ \sum_{\tau=t}^{t+T-1} \bbbe \Big[ \frac{\Delta V(\mathbf{X}(\tau)) - \Delta V_{\parallel}(\mathbf{X}(\tau))}{2\|\mathbf{q}_{\perp}(\tau)\|} \Big| \mathbf{X}(\tau) \Big] \bigg] 
\label{drift_q_perp_1}
\end{align}

We now derive bounds for $\Delta V$ and $\Delta V_{\parallel}$:
\begin{align*}
\bbbe \Big[ \Delta V(\mathbf{X}(\tau)) \Big| \mathbf{X}(\tau) \Big]
=& \bbbe \Big[ \|\mathbf{q}(\tau+1)\|^2 - \|\mathbf{q}(\tau)\|^2 \Big| \mathbf{X}(\tau) \Big] \\
=& \bbbe \Big[ \|\mathbf{q}(\tau) + \mathbf{a}(\tau) - \mathbf{s}(\tau)\mathds{1}_{\{r(\tau) = 0\}} \|^2 + \|\mathbf{u}(\tau) \|^2 + 2\Big\langle \mathbf{q}(\tau+1) - \mathbf{u}(\tau), \mathbf{u}(\tau) \Big\rangle - \|\mathbf{q}(\tau)\|^2 \Big| \mathbf{X}(\tau)  \Big] \\
\leq& \bbbe \Big[ \|\mathbf{q}(\tau) + \mathbf{a}(\tau) - \mathbf{s}(\tau)\mathds{1}_{\{r(\tau) = 0\}} \|^2 - \|\mathbf{q}(\tau)\|^2  \Big| \mathbf{X}(\tau)  \Big] \\
=& \sum_{i,j} \bbbe \Big[ a_{ij}^2(\tau) + s_{ij}(\tau)\mathds{1}_{\{r(\tau) = 0\}} - 2 a_{ij}(\tau) s_{ij}(\tau) \mathds{1}_{\{r(\tau) = 0\}} \Big| \mathbf{X}(\tau)  \Big] + \bbbe \Big[ 2 \Big\langle \mathbf{q}(\tau), \boldsymbol{\lambda} -  \mathbf{s}(\tau)\mathds{1}_{\{r(\tau) = 0\}} \Big\rangle \Big| \mathbf{X}(\tau)  \Big]  \\
\stackrel{(a)}{\leq}& \sum_{ij} (\lambda_{ij}^2 + \sigma_{ij}^2) + n +  2 \Big\langle \mathbf{q}(\tau), (1-\epsilon)\boldsymbol{\nu} -  \mathbf{s}(\tau)  \Big\rangle + 2 \Big\langle \mathbf{q}(\tau), \mathbf{s}(\tau)  \Big\rangle \mathds{1}_{\{r(\tau) > 0\}}  \\
=& \|\boldsymbol{\lambda}\|^2 + \|\boldsymbol{\sigma}\|^2 + n - 2\epsilon \Big\langle \mathbf{q}(\tau), \boldsymbol{\nu} \Big\rangle + 2 \Big\langle \mathbf{q}(\tau), \boldsymbol{\nu} - \mathbf{s}(\tau) \Big\rangle  + 2 \Big\langle \mathbf{q}(\tau), \mathbf{s}(\tau)  \Big\rangle \mathds{1}_{\{r(\tau) > 0\}} 
\end{align*}
where $(a)$ follows from $\bbbe[ a_{ij}^2] = \lambda_{ij}^2 + \sigma_{ij}^2$, $a_{ij}(t)s_{ij}(t) \geq 0$ for all $i,j$, and $\sum_{ij} s_{ij}(t) = 1$ for all $t$.

Suppose $g$ is the sublinear hysteresis function for the Adaptive MaxWeight, then by the sublinearity, there exists a constant $K_{\theta}$ such that $g(x) < \frac{\theta}{\alpha} x$ for any $x > K_{\theta}$, where $\alpha = \frac{8\|\boldsymbol{\nu}\|}{\nu_{\min}}$. Hence by the definition of the Adaptive MaxWeight, we have for any $\mathbf{X}(\tau)$ such that $\langle \mathbf{q}(\tau), \mathbf{s}^*(\tau) \rangle > K_{\theta}$:
\begin{align*}
\Big\langle \mathbf{q}(\tau),  \boldsymbol{\nu}-\mathbf{s}(\tau)  \Big\rangle 
=& \Big\langle \mathbf{q}(\tau),  \boldsymbol{\nu}-\mathbf{s}^*(\tau)  \Big\rangle 
+ \Big\langle \mathbf{q}(\tau),  \mathbf{s}^*(\tau)-\mathbf{s}(\tau)  \Big\rangle  \\
\leq& \Big\langle \mathbf{q}(\tau), \boldsymbol{\nu}-\mathbf{s}^*(\tau) \Big\rangle + g\left( \Big\langle \mathbf{q}(\tau), \mathbf{s}^*(\tau) \Big\rangle \right)  \\
\leq& \Big\langle \mathbf{q}(\tau), \boldsymbol{\nu}-\mathbf{s}^*(\tau) \Big\rangle + \frac{\theta}{\alpha} \Big\langle \mathbf{q}(\tau), \mathbf{s}^*(\tau) \Big\rangle  \\
=& \Big(1- \frac{\theta}{\alpha} \Big)  \Big\langle \mathbf{q}(\tau), \boldsymbol{\nu}-\mathbf{s}^*(\tau) \Big\rangle + \frac{\theta}{\alpha} \Big\langle \mathbf{q}(\tau), \boldsymbol{\nu} \Big\rangle
\end{align*}

From~\cite[Claim 2]{MW_scaling}, we have $\Big\langle \mathbf{q}(\tau), \boldsymbol{\nu}-\mathbf{s}^*(\tau) \Big\rangle \leq -\nu_{\min} \|\mathbf{q}_{\perp}(\tau)\|$. Therefore,
\begin{align}
\bbbe \Big[ \Delta V(\mathbf{X}(\tau)) \Big| \mathbf{X}(\tau) \Big]
\leq& \|\boldsymbol{\lambda}\|^2 + \|\boldsymbol{\sigma}\|^2 + n - 2\epsilon \Big\langle \mathbf{q}(\tau), \boldsymbol{\nu} \Big\rangle -2 \Big(1 - \frac{\theta}{\alpha}\Big) \nu_{\min}\|\mathbf{q}_{\perp}(\tau)\|  \nonumber \\
&+ 2 \frac{\theta}{\alpha} \Big\langle \mathbf{q}(\tau), \boldsymbol{\nu} \Big\rangle   + 2 \Big\langle \mathbf{q}(\tau), \mathbf{s}(\tau)  \Big\rangle \mathds{1}_{\{r(\tau) > 0\}} 
\label{drift_V}
\end{align}

For $\Delta V_{\parallel}$, we have
\begin{align}
&\bbbe \Big[ \Delta V_{\parallel}(\mathbf{X}(\tau)) \Big| \mathbf{X}(\tau) \Big] = \bbbe \Big[ \|\mathbf{q}_{\parallel}(\tau+1)\|^2 - \|\mathbf{q}_{\parallel}(\tau)\|^2 \Big| \mathbf{X}(\tau) \Big]  \nonumber \\
=& \bbbe \Big[ \Big\langle \mathbf{q}_{\parallel}(\tau+1) + \mathbf{q}_{\parallel}(\tau), \mathbf{q}_{\parallel}(\tau+1) - \mathbf{q}_{\parallel}(\tau) \Big\rangle  \Big| \mathbf{X}(\tau) \Big]  \nonumber \\
=& \bbbe \Big[ \|\mathbf{q}_{\parallel}(\tau+1) - \mathbf{q}_{\parallel}(\tau)\|^2 + 2\Big\langle  \mathbf{q}_{\parallel}(\tau), \mathbf{q}_{\parallel}(\tau+1) - \mathbf{q}_{\parallel}(\tau) \Big\rangle \Big| \mathbf{X}(\tau) \Big]  \nonumber \\
\geq&  2\bbbe \Big[ \Big\langle  \mathbf{q}_{\parallel}(\tau), \mathbf{q}_{\parallel}(\tau+1) - \mathbf{q}_{\parallel}(\tau) \Big\rangle \Big| \mathbf{X}(\tau) \Big]    \nonumber \\
=&  2\bbbe \Big[ \Big\langle  \mathbf{q}_{\parallel}(\tau), \mathbf{q}(\tau+1) - \mathbf{q}(\tau) \Big\rangle - \Big\langle  \mathbf{q}_{\parallel}(\tau), \mathbf{q}_{\perp}(\tau+1) - \mathbf{q}_{\perp}(\tau) \Big\rangle \Big| \mathbf{X}(\tau) \Big] \nonumber \\
\stackrel{(b)}{\geq}&  2\bbbe \Big[ \Big\langle  \mathbf{q}_{\parallel}(\tau), \mathbf{a}(\tau) - \mathbf{s}(\tau) \mathds{1}_{\{r(\tau) = 0\}} + \mathbf{u}(\tau) \Big\rangle \Big| \mathbf{X}(\tau) \Big]  \nonumber \\
\geq&  2\Big\langle \mathbf{q}_{\parallel}(\tau), \boldsymbol{\lambda} \Big\rangle -2 \Big\langle \mathbf{q}_{\parallel}(\tau), \mathbf{s}(\tau) \mathds{1}_{\{r(\tau) = 0\}} \Big\rangle  \nonumber \\
=& -2\epsilon \Big\langle \mathbf{q}_{\parallel}(\tau), \boldsymbol{\nu} \Big\rangle + 2 \Big\langle \mathbf{q}_{\parallel}(\tau), \boldsymbol{\nu} - \mathbf{s}(\tau) \Big\rangle + 2 \Big\langle \mathbf{q}_{\parallel}(\tau), \mathbf{s}(\tau)  \Big\rangle \mathds{1}_{\{r(\tau) > 0\}}  \nonumber \\
=& -2\epsilon \Big\langle \mathbf{q}_{\parallel}(\tau), \boldsymbol{\nu} \Big\rangle + 2 \Big\langle \mathbf{q}_{\parallel}(\tau), \mathbf{s}(\tau)  \Big\rangle \mathds{1}_{\{r(\tau) > 0\}} 
\label{drift_V_parallel}
\end{align}
For $(b)$, we use the following properties of the projection onto cone $\mathcal{K}$. For $\mathbf{q} \in \bbbr^{n^2}$, $\langle \mathbf{q}_{\parallel}, \mathbf{q}_{\perp} \rangle = 0$, and $\mathbf{q}_{\perp} \in \mathcal{K}^{\circ}$. Therefore, $\langle \mathbf{q}_{\parallel}(t), \mathbf{q}_{\perp}(t) \rangle = 0$, and $\langle \mathbf{q}_{\parallel}(t), \mathbf{q}_{\perp}(t+1) \rangle \leq 0$.

Apply (\ref{drift_V}) and (\ref{drift_V_parallel}) into (\ref{drift_q_perp_1}), we obtain
\begin{align}
\bbbe_{\mathbf{X}} \Big[ \|\mathbf{q}_{\perp}(t+T)\| - \|\mathbf{q}_{\perp}(t)\|  \Big]
\leq& \bbbe_{\mathbf{X}} \bigg[ \sum_{\tau=t}^{t+T-1}  \Big( \frac{\|\boldsymbol{\lambda}\|^2 + \|\boldsymbol{\sigma}\|^2 + n}{2\|\mathbf{q}_{\perp}(\tau)\|}  - \epsilon \Big\langle \frac{\mathbf{q}_{\perp}(\tau)}{\|\mathbf{q}_{\perp}(\tau)\|}, \boldsymbol{\nu} \Big\rangle  -  \big(1-\frac{\theta}{\alpha}\big) \nu_{\min}  \nonumber \\
&+ \frac{\theta  \big\langle \mathbf{q}(\tau), \boldsymbol{\nu} \big\rangle}{\alpha \|\mathbf{q}_{\perp}(\tau)\|} + \Big\langle \frac{\mathbf{q}_{\perp}(\tau)}{\|\mathbf{q}_{\perp}(\tau)\|}, \mathbf{s}(\tau) \Big\rangle \mathds{1}_{\{r(\tau) > 0\}}  \Big) \bigg]  \nonumber \\
\leq& \bbbe_{\mathbf{X}} \bigg[ T \bigg( \frac{\|\boldsymbol{\lambda}\|^2 + \|\boldsymbol{\sigma}\|^2 + n}{\min\limits_{\tau \in [t, t+T]}2\|\mathbf{q}_{\perp}(\tau)\|} + \epsilon \|\boldsymbol{\nu}\| -  \big(1-\theta \big) \nu_{\min} + \frac{1 + \theta}{\alpha} \|\boldsymbol{\nu}\| \bigg) +  \sqrt{n}\sum_{\tau=t}^{t+T-1}\mathds{1}_{\{r(\tau) > 0\}}  \bigg] 
\label{drift_q_perp_2}
\end{align}
where we have used the fact that $\|\mathbf{q}_{\perp}\| \geq \theta \|\mathbf{q}_{\parallel}\|$ implies $\theta \|\mathbf{q}\| \leq \theta (\|\mathbf{q}_{\parallel}\| + \|\mathbf{q}_{\perp}\|) \leq (1+\theta) \|\mathbf{q}_{\perp}$, and $\|\mathbf{s}(\tau)\| \leq \sqrt{n}$ for any schedule $\mathbf{s}(\tau) \in \mathcal{S}$.

On the other hand, the drift of $\|\mathbf{q}_{\parallel}\|$ could be obtained following (\ref{drift_V_parallel}):
\begin{align}
\bbbe_{\mathbf{X}} \Big[ \|\mathbf{q}_{\parallel}(t+T)\| - \|\mathbf{q}_{\parallel}(t)\|  \Big]
=&  \bbbe_{\mathbf{X}} \bigg[ \sum_{\tau=t}^{t+T-1} \bbbe \Big[ \|\mathbf{q}_{\parallel}(\tau+1)\| - \|\mathbf{q}_{\parallel}(\tau)\| \Big| \mathbf{X}(\tau) \Big]  \bigg] \nonumber \\
\geq& \bbbe_{\mathbf{X}} \bigg[ \sum_{\tau=t}^{t+T-1} \bbbe \Big[ \frac{\|\mathbf{q}_{\parallel}(\tau+1)\|^2 - \|\mathbf{q}_{\parallel}(\tau)\|^2}{\|\mathbf{q}_{\parallel}(\tau+1)\|+\|\mathbf{q}_{\parallel}(\tau)\|} \Big| \mathbf{X}(\tau) \Big]  \bigg] \nonumber \\
\geq& \bbbe_{\mathbf{X}} \bigg[ \sum_{\tau=t}^{t+T-1} \bbbe \Big[ \frac{-2\epsilon \big\langle \mathbf{q}_{\parallel}(\tau), \boldsymbol{\nu} \big\rangle}{\|\mathbf{q}_{\parallel}(\tau+1)\|+\|\mathbf{q}_{\parallel}(\tau)\|} \Big| \mathbf{X}(\tau) \Big]  \bigg] \nonumber \\
\geq& \bbbe_{\mathbf{X}} \bigg[ \sum_{\tau=t}^{t+T-1} \frac{-2\epsilon \big\langle \mathbf{q}_{\parallel}(\tau), \boldsymbol{\nu} \big\rangle}{\|\mathbf{q}_{\parallel}(\tau)\|}   \bigg]  \geq -2T \epsilon \|\boldsymbol{\nu}\|
\label{drift_q_parallel}
\end{align}
where the last inequality follows the Cauchy-Schwartz inequality.

Now apply (\ref{drift_q_perp_2}) and (\ref{drift_q_parallel}) into (\ref{drift_Z}), we obtain
\begin{align}
\bbbe_{\mathbf{X}} \Big[ \Delta^T Z(\mathbf{X})  \Big]
\leq& \bbbe_{\mathbf{X}} \bigg[ T \bigg( \frac{\|\boldsymbol{\lambda}\|^2 + \|\boldsymbol{\sigma}\|^2 + n}{\min\limits_{\tau \in [t, t+T]}2\|\mathbf{q}_{\perp}(\tau)\|} + (1+2\theta)\epsilon \|\boldsymbol{\nu}\| -  \big(1-\theta \big) \nu_{\min} + \frac{1 + \theta}{\alpha} \|\boldsymbol{\nu}\| \bigg) +  \sqrt{n}\sum_{\tau=t}^{t+T-1}\mathds{1}_{\{r(\tau) > 0\}}  \bigg] 
\end{align}

From~\cite[Lemma 1]{AMW_arxiv}, we know that for any fixed $T>0$, if $W^*(t) > g^{-1}\Big( nT(a_{\max}+1) \Big) + nT$, then at most one reconfiguration could occur within $[t, t+T]$, which gives $\sum_{\tau=t}^{t+T-1} \mathds{1}_{r(\tau) > 0} \leq \Delta_r$.

Select $T = \frac{8 \sqrt{n} \Delta_r}{\nu_{\min}}$, then set $D = \frac{3}{2}na_{\max}T = \frac{12 n^{3/2} a_{\max} \Delta_r}{\nu_{\min}}$. Then $\forall \mathbf{X}: Z(\mathbf{X}) > \kappa = \Big\{ D, na_{\max}T + \frac{4(\|\boldsymbol{\lambda}\|^2 + \|\boldsymbol{\sigma}\|^2 + n)}{\nu_{\min}} , nK_{\theta},$ $ng^{-1}\Big( nT(a_{\max}+1) \Big) + n^2T \Big\}$, and $\forall \epsilon : 0 < \epsilon \leq \frac{\nu_{\min}}{16\|\boldsymbol{\nu}\|}$, we have $\bbbe_{\mathbf{X}} \Big[ \Delta^T Z(\mathbf{X}) \Big] \leq -\frac{(1-\theta)\nu_{\min}}{4} \leq -\frac{\nu_{\min}}{8}$.

Hence by Lemma~\ref{lemma:Foster}, we have $\forall \epsilon: 0 < \epsilon \leq \frac{\nu_{\min}}{16\|\boldsymbol{\nu}\|}$:
\begin{align*}
\bbbe \left[ \| \bar{\mathbf{q}}_{\perp} \| - \theta \| \bar{\mathbf{q}}_{\parallel} \| \right]
\leq \bbbe \left[ Z(\bar{\mathbf{e}}) \right]
\leq \kappa + \frac{16 D^2}{\nu_{\min}}.
\end{align*}
Let $M_{\theta} = \kappa + \frac{16 D^2}{\nu_{\min}}$, we then have the result.

\end{proof}

\end{document}